\documentclass[]{amsart}
\usepackage{amssymb}
\usepackage{latexsym}
\input epsf
\newtheorem{prop}{Proposition}
\newtheorem{thm}{Theorem}

\def\P{{\mathbb P}}

\def\R{{\mathbb R}}
\def\C{{\mathbb C}}

\def\H{{\mathcal H}}
\def\l{{\lambda}}
\def\m{{\mu}}

\def\L{{\Lambda}}
\def\a{{\alpha}}
\def\b{{\beta}}
\def\g{{\gamma}}
\def\s{{\sigma}}
\def\om{{\omega}}
\def\Om{{\Omega}}
\def\endproof{\hfill $\Box$}

\newcommand{\dis}{\displaystyle}
\newcommand{\partialr}{\partial}

\newcommand{\gequ}{\geqslant}
\newcommand{\lequ}{\leqslant}
\newcommand{\ra}{\rightarrow}
\newcommand{\skipline}{\vspace{\baselineskip}}

\newcommand{\noin}{\noindent}
\newcommand{\wt}{\widetilde}

\newcommand{\lie}[1]{{\mathfrak #1}}
\newcommand{\refg}[1]{{(\ref{#1})}}

\newcommand{\cdott}{*}

\begin{document}

\title{The Hodge Star Operator on Schubert Forms}
\author{Klaus K\"{u}nnemann and Harry Tamvakis}
\date{submitted July 14, 2000, accepted January 23, 2001, to appear in {\it Topology}}
\address{NWF-I Mathematik, Universit\"at Regensburg, 93040 Regensburg,
Germany}
\email{klaus.kuennemann@mathematik.uni-regensburg.de}
\address{Department of Mathematics, University of Pennsylvania,
209 South 33rd Street, 
Philadelphia, PA 19104-6395, USA}
\email{harryt@math.upenn.edu}
\begin{abstract}
Let $X=G/P$ be a homogeneous space of a complex semisimple
Lie group $G$ equipped with a hermitian metric. 
We study the action of the Hodge star operator
on the space of harmonic differential forms on $X$. We obtain
explicit combinatorial formulas for this action when $X$
is an irreducible hermitian symmetric space of compact type.
\end{abstract}

\maketitle 

\setcounter{section}{-1}

\section{Introduction}
\noindent
Let us recall the definition of the Hodge $*$-operator. If $V$ is
an $n$-dimensional Euclidean vector space, choose an orthonormal basis 
$e_1,\ldots,e_n$ of $V$ and define the star operator
$*:\wedge^kV\ra\wedge^{n-k}V$ by
\[
\dis
*\,(e_{\s(1)}\wedge\ldots\wedge e_{\s(k)})=\mathrm{sgn}(\s)\,
e_{\s(k+1)}\wedge\ldots\wedge e_{\s(n)}
\]
for any permutation $\sigma$ of the indices $(1,\ldots,n)$.
The operator $*$ depends only on the inner product
structure of $V$ and the orientation determined by the basis
$e_1,\ldots,e_n$. 
If $V$ is the real vector space underlying a hermitian space then $*$
is defined  using the natural choice of orientation coming from the 
complex structure.

Let $X$ be a hermitian complex manifold of complex dimension $d$, and 
 $A^{p,q}(X)$ the space of complex valued smooth differential forms of type
 $(p,q)$ on $X$. The star
operator taken pointwise gives a complex linear isomorphism
$*:A^{p,q}(X)\ra A^{d-q,d-p}(X)$ such that $**=(-1)^{p+q}$. Since 
$*$ commutes with the Laplacian $\Delta$, it induces an
isomorphism $\H^k(X)\ra\H^{2d-k}(X)$ between spaces of harmonic forms 
on $X$. This in turn gives a map on cohomology groups
$*:H^k(X,\C)\ra H^{2d-k}(X,\C)$ which depends on the metric structure
of $X$.

Our main goal is to compute the action of $*$ explicitly when $X$ is an
irreducible hermitian symmetric space of compact type, equipped with
an invariant K\"{a}hler metric. These spaces have been classified by
\'{E}.\ Cartan; there are four infinite families and two `exceptional'
cases. We will describe our result here in the case of the
Grassmannian 
\[
\dis
G(m,n)=U(m+n)/(U(m)\times U(n))
\]
of complex $m$-dimensional linear subspaces of $\C^{m+n}$.

\begin{figure}[htpb]
\epsfxsize 35mm
\center{\mbox{\epsfbox{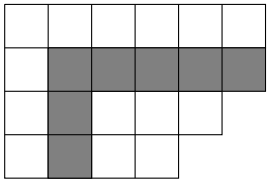}}} 
\caption{A hook in the diagram of $(6,6,5,4)$}
\label{hook}
\end{figure}

Recall that a partition $\l=(\l_1\gequ\l_2\gequ\cdots\gequ\l_k)$ 
is identified with its Young diagram of boxes; the {\em weight} $|\l|$ 
is the number of boxes in $\l$. 
%Thus $\l$ is a partition of the
%number $k=|\l|$, denoted $\l\vdash k$.
Given a diagram $\l$ and a box $x\in\l$, the {\em hook} $H_x$ is the
set of all boxes directly to the right and below $x$, including
$x$ itself (see Figure \ref{hook}).
The number of boxes in $H_x$ is the {\em hook length} $h_x$. 
We let
\[
\dis
h_{\l}:=\prod_{x\in\l}h_x
\]
denote the product of the hook lengths in $\l$. It is known that
$N_{\l}=|\l|!/h_{\l}$ is the dimension of the irreducible
representation of the symmetric group $S_{|\l|}$ corresponding to
$\l$. The integer
$N_{\l}$ also counts the number of {\em standard Young tableaux}
on $\l$, that is, the number of different ways to fill the boxes in 
$\l$ with the numbers $1,2,\ldots, |\l|$ so that the entries are 
strictly increasing along rows and columns. This fact is due to
Frame, Robinson, and Thrall \cite{FRT}.

\begin{figure}[htpb]
\epsfxsize 35mm
\center{\mbox{\epsfbox{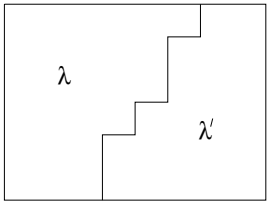}}} 
\caption{The diagrams for Poincar\'e dual forms on $G(6,8)$}
\label{duals}
\end{figure}

Partitions parametrize the harmonic forms corresponding to the
Schubert classes, which are the natural geometric
basis for the cohomology ring of $G=G(m,n)$. 
For each partition $\l$ whose diagram is 
contained in the $m\times n$ rectangle $(n^m)$, there is a 
harmonic form $\Omega_{\l}$ 
of type $(|\l|,|\l|)$ which is dual to the class
of the codimension $|\l|$ Schubert variety $X_{\l}$ in $G$. 
The Poincar\'{e}
dual of $\Omega_{\l}$ (i.e. the dual form with respect to the
Poincar\'{e} pairing
$(\phi,\psi)\longmapsto \int_G\phi\wedge\psi$)
corresponds to the diagram $\l'$
which, when inverted, is the complement of $\l$ in $(n^m)$ (see Figure
\ref{duals}).

Normalize a given invariant hermitian metric on $G$ so that its fundamental 
form is the Schubert form $\Om_1$. 
We can now state our result for the action of the star operator:
\begin{equation}
\label{introeq}
*\,\Omega_{\l}=\frac{h_{\l}}{h_{\l'}}\,\Omega_{\l'}.
\end{equation}
There are similar results for infinite families of different type:
the even orthogonal and Lagrangian Grassmannians. In 
these cases combinatorialists have identified the correct notions of
`hook' and `hook length', and our formula for $*$ is a direct
analogue of (\ref{introeq}). We also compute the action of $*$ for
quadric hypersurfaces and the exceptional cases.

More generally, our calculations are valid for any K\"ahler manifold
whose cohomology ring coincides with that of a hermitian symmetric space.
For example we compute the action of $*$ on the harmonic forms 
(with respect to any K\"ahler metric) for the odd orthogonal 
Grassmannians $SO(2n+1)/U(n)$. 

The motivation for this work came from Arakelov geometry. A combinatorial
understanding of the Lefschetz theory on homogeneous spaces is useful
in the study of the corresponding objects over the ring of integers.
For such arithmetic varieties, Gillet and Soul\'e have formulated
analogues of Grothendieck's standard conjectures on 
algebraic cycles (see \cite[Section 5.3]{So}). Our calculation of
$*$ has been used by Kresch and the second author \cite{KT} to
verify these conjectures for the arithmetic Grassmannian $G(2,n)$. 

Let us briefly outline the contents of the paper.
In the first section we state our main theorem about the action of the
Hodge star operator on Schubert forms on irreducible compact hermitian
symmetric spaces.
In the next section, we show that the star operator 
on these spaces
maps a Schubert form to a non-zero multiple of the Poincar\'e dual form.
Our proof is based on Kostant's construction of a basis of the space of 
harmonic forms which is recalled here.
We prove our main theorem in section three.
In the fourth section, we use our calculation of the star operator to give
an explicit formula for the adjoint of the Lefschetz operator.
Here we recover some results of Proctor.
The example in the final section shows that for complete flag
varieties, the star operator
no longer maps a Schubert form to a multiple of 
the Poincar\'e dual form.

Work on this paper started while the authors were visiting the Isaac
Newton Institute in Cambridge.
We greatfully acknowledge the hospitality of the Institute.
It is a pleasure to thank the Deutsche Forschungsgemeinschaft (first
author) and the National Science Foundation (second author) for support
during the preparation of the paper. We also thank
Wolfgang Ziller for several valuable discussions and Anton Deitmar
for his comments on a first version of this text.

%
%\bigskip
%\bigskip
%
%\begin{center}
%\sc Contents
%\end{center}
%
%\medskip
%
%
%\contentsline {section}{\tocsection {}{0}{Introduction}}{1}
%\contentsline {section}{\tocsection {}{1}{Statement of the Main Theorem}}{3}
%\contentsline {section}{\tocsection {}{2}{Schubert Forms on 
%Homogeneous Spaces}}{6}
%\contentsline {section}{\tocsection {}{3}{Proof of the Main Theorem}}{8}
%\contentsline {section}{\tocsection {}{4}{Formulas for the Adjoint of the 
%Lefschetz Operator}}{10}
%\contentsline {section}{\tocsection {}{5}{An Example}}{13}
%\contentsline {section}{\tocsection {}{}{References}}{14}
%
%
%\bigskip
%

\section{Statement of the Main Theorem}
\label{mt}
We begin with some more notation from combinatorics: a partition
$\l=(\l_i)_{i\gequ 1}$ is {\em strict} if its
parts $\l_i$ are distinct; the number of non-zero parts is the
{\em length} of $\l$, denoted $\ell(\l)$. Define
$\alpha(\l)=|\l|-\ell(\l)$ to be the number of boxes in $\l$ that
are not in the first column.

\begin{figure}[htpb]
\epsfxsize 35mm
\center{\mbox{\epsfbox{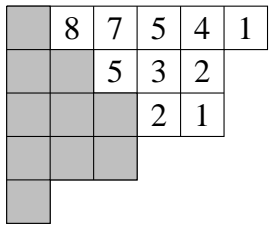}}} 
\caption{$D(\l)$, $S(\l)$ and hook lengths for $\l=(5,3,2)$}
\label{double}
\end{figure}

For a strict partition $\l=(\l_1>\l_2>\cdots >\l_m>0)$, the
{\em shifted diagram} $S(\l)$ is obtained from the Young diagram
of $\l$ by moving the $i$th row $(i-1)$ squares to the right, for 
each $i>1$. The 
{\em double diagram} $D(\l)$ consists of $S(\l)$ dovetailed into its
reflection in the main diagonal $\{(i,i): i>0\}$; in 
Frobenius notation, we have
$D(\l)=(\l_1,\ldots,\l_m\,|\, \l_1-1,\ldots,\l_m-1)$
(this is illustrated in Figure \ref{double}; see also
\cite[Section I.1]{M} for Frobenius notation).
For each box $x$ in $S(\l)$, the hook length $h_x$ is
defined to be the hook length at $x$ in the double diagram $D(\l)$.
Figure \ref{double} displays these hook lengths for $\l=(5,3,2)$.

We let
\[
\dis
g_{\l}:=\prod_{x\in S(\l)}h_x,
\]
the product over all boxes $x$ of the shifted diagram of $\l$.
We remark that $|\l|!/g_{\l}$ counts the number of standard shifted
tableaux of shape $\l$; it also occurs in the degree formula for
the corresponding projective representation of the symmetric group.
See for instance \cite[p.\ 187 and Theorem 10.7]{HH} for definitions
and details.

The irreducible compact hermitian symmetric spaces have been
classified by \'E.\ Cartan \cite{C}; there are four infinite families
and two exceptional cases (see \cite{Wo} for a modern treatment, in
particular Corollary 8.11.5). 
We will recall this list; each such space is of the form $G/P=K/V$, 
with notation as in Section \ref{sosc} (we will use the compact
presentation $K/V$). In each case we have a natural $K$-invariant hermitian
metric, unique up to positive scalar.
The $K$-invariant differential forms coincide with the
harmonic forms for this metric, and are all of
$(p,p)$ type for some $p$.

The Schubert cycles form a basis for the integral homology ring;
their duals in cohomology are represented by unique harmonic forms,
called {\em Schubert forms}. The Schubert forms are parametrized by
the set $W^1$ defined in (\ref{w1}) below, and for the first three
infinite families that follow, this parameter space can be 
realized using integer partitions. We refer to \cite{BGG}, \cite{CN},
\cite{Hiller}, \cite[Section 12]{Pro2} for more information.

\medskip
\noin
(i) The Grassmannian $G(m,n)=U(m+n)/(U(m)\times U(n))$ of 
$m$-dimensional linear subspaces of $\C^{m+n}$, with
$\dim_{\C}G(m,n)=mn$. We have a Schubert
form $\Omega_{\l}$ of type $(|\l|,|\l|)$ 
for each partition $\l$ whose Young diagram is
contained in the $m\times n$ rectangle $(n^m)=(n,\ldots,n)$.
The Poincar\'e dual form corresponds to the diagram $\l'$
described in the introduction.

\medskip
\noin
(ii) The even orthogonal Grassmannian $OG(n,2n)=SO(2n)/U(n)$ 
(spinor variety) parametrizing
maximal isotropic subspaces of $\C^{2n}$ equipped with a nondegenerate
symmetric form, with $\dim_{\C}OG(n,2n)=\binom{n}{2}$.
There is
a Schubert form $\Phi_{\l}$ of type $(|\l|,|\l|)$ 
for each strict partition $\l$
whose diagram is contained in the triangular partition
$\rho_{n-1}:=(n-1,n-2,\ldots,1)$.
The Poincar\'{e}
dual of $\Phi_{\l}$ corresponds to the strict partition $\l'$ whose 
parts complement the parts of $\l$ in the set $\{1,\ldots,n-1\}$.

\medskip
\noin
(iii) The Lagrangian
Grassmannian $LG(n,2n)=Sp(2n)/U(n)$ parametrizing Lagrangian
subspaces of $\C^{2n}$ equipped with a symplectic form, with
$\dim_{\C}LG(n,2n)=\binom{n+1}{2}$. 
Here we have a Schubert form of type $(|\l|,|\l|)$,
denoted $\Psi_{\l}$,
for each strict partition $\l$ whose diagram is contained in 
$\rho_n=(n,\ldots,1)$.
The Poincar\'{e}
dual of $\Psi_{\l}$ corresponds to the strict partition $\l'$ whose
parts complement the parts of $\l$ in the set $\{1,\ldots,n\}$.

\medskip
\noin
(iv) The complex quadric $Q(n)=SO(n+2)/(SO(n)\times SO(2))$, isomorphic
to a smooth quadric hypersurface in $\P^{n+1}(\C)$, of dimension $n$.
Let $\omega$ denote the K\"ahler form which is dual to the
class of a hyperplane. 
Our reference for the cohomology ring of $Q(n)$
is \cite[Section 2]{EG} (working in the context of Chow rings).
There are two cases depending on the parity of $n$:

\medskip
\noin
{\bf --}\ If $n=2k-1$ is odd then there is one Schubert form $e \in
\H^{2k}(X)$ and $\om^k=2e$. The complete list of Schubert forms is
\[
\dis
1,\om,\om^2,\ldots,\om^{k-1},e,\om e,\ldots, \om^{k-1}e.
\]
The Poincar\'{e} dual of $\om^i$ is $\om^{k-1-i}e$, for 
$0\lequ i\lequ k-1$.

\medskip
\noin
{\bf --}\  If $n=2k$ is even then there are two distinct
Schubert forms $e_0,e_1 \in
\H^{2k}(X)$ which correspond to the two rulings of the
quadric (in homology). Furthermore $\om^k=e_0+e_1$ and the
complete list of Schubert forms is
\[
\dis
1,\om,\om^2,\ldots,\om^{k-1},e_0,e_1,\om e_0=\om e_1,\om^2 e_0,
\ldots, \om^ke_0.
\]
The Poincar\'{e} dual of $\om^i$ is $\om^{k-i}e_0$, for 
$0\lequ i\lequ k-1$, while the dual of $e_j$ is $e_{j+k}$, where
the indices are taken mod 2.

\medskip
\noin
(v) The `exceptional' space $E_6/(SO(10)\cdot SO(2))$, of complex
dimension 16.

\medskip
\noin
(vi) The `exceptional' space $E_7/(E_6\cdot SO(2))$, of complex 
dimension 27.

\medskip

Suppose $X$ is a compact K\"ahler manifold whose cohomology ring is
isomorphic to any occuring in the preceding examples. 
We will see in Section \ref{comb} that the action of the star operator
on $H^*(X,\C)$ is the same as if $X$ were a hermitian symmetric space.
If we look among the homogeneous spaces $G/P$ considered in 
Section \ref{sosc} we find one example with this property which is not
itself hermitian symmetric:

\medskip
\noin
(ii$'$) The odd orthogonal Grassmannian $OG(n-1,2n-1)=SO(2n-1)/U(n-1)$ 
parametrizing
maximal isotropic subspaces of $\C^{2n-1}$ equipped with a nondegenerate
symmetric form, whose cohomology ring coincides with that of 
$OG(n,2n)$. Choose any K\"ahler metric on
$OG(n-1,2n-1)$. By abuse of notation
we use $\Phi_{\l}$ to denote the harmonic Schubert form corresponding
to the strict partition $\l\subset\rho_{n-1}$; in this way the
statement of the next theorem will include this example.

\skipline

Let $\omega$ denote the fundamental form of a hermitian metric $h$ on $X$
which is given in any local holomorphic coordinate chart $(z_i)$ as
\[
\omega=\frac{i}{2}\,\sum_{i,j}\,h\,\biggl(\frac{\partial}{\partial z_i},
\frac{\partial}{\partial z_j}\biggr)\, dz_i\wedge d\overline{z}_j.
\]
In each example we normalize the hermitian metric so that its
fundamental form coincides with the unique Schubert form of type 
$(1,1)$. 
We can now state our main result computing the Hodge star operator in
examples (i)-(iv). The two exceptional cases will be discussed in 
Section \ref{proofsec}.

\begin{thm}
\label{mainthm}
The action of the Hodge star operator on the Schubert forms in examples
(i)--(iii) is as follows:
\[
\dis
*\,\Omega_{\l}=\frac{h_{\l}}{h_{\l'}}\,\Omega_{\l'},\ \ \
*\,\Phi_{\l}=\frac{g_{\l}}{g_{\l'}}\,\Phi_{\l'},\ \ \ 
*\,\Psi_{\l}=2^{\alpha(\l')-\alpha(\l)}
\,\frac{g_{\l}}{g_{\l'}}\,\Psi_{\l'}.
\]
In example (iv) if $n=2k-1$ is odd we have
\[
\dis
*\,\om^i=\frac{2\cdot i!}{(n-i)!}\,\om^{k-1-i}e, \ \ \ 0\lequ i\lequ k-1
\]
and if $n=2k$ is even then
\[
\dis
*\,\om^i=\frac{2\cdot i!}{(n-i)!}\,\om^{k-i}e_0, \ \ \ 0\lequ i\lequ k-1
\ \ \ \ and \ \ \ *e_j=e_{j+k},
\]
while the remaining terms are determined by the relation $**=1$.

\end{thm}

\section{Schubert Forms on Homogeneous Spaces}
\label{sosc}
We recall some fundamental facts from Kostant's seminal
papers \cite{1}, \cite{2}.
We derive from Kostant's results that a Schubert form 
on an irreducible compact hermitian symmetric space is mapped by the Hodge
star operator to a multiple of the dual Schubert form.
We consider the following situation:
 $\lie{g}$  is a complex semi-simple Lie algebra,
 $\lie{p}$  a parabolic subalgebra of $\lie{g}$ which contains a fixed
Borel subalgebra $\lie{b}$,
 $\lie{n}$ the maximal nilpotent ideal of $\lie{p}$, and
 $\lie{k}$ a fixed compact real form of $\lie{g}$.
The choice of $\lie{k}$ determines a Cartan involution $x\mapsto x^{\theta}$  on
$\lie{g}$ defined as $(u+iv)^{\theta}=u-iv$ for $u,v\in \lie{k}$.
  % \cite[3.3]{2}.
For any subspace $\lie{s}$ of $\lie{g}$, we set 
$\lie{s}^{\theta}=\{x^{\theta}\, |\, x\in \lie{s}\}$.
Let $\lie{g}_1=\lie{p}\cap \lie{p}^{\theta}\subset\lie{g}$.
We have $\lie{p}=\lie{g}_1+\lie{n}$ and
 $\lie{g}=\lie{n}+\lie{g}_1+\lie{n}^{\theta}$.
Let $\lie{r}=\lie{n}+\lie{n}^{\theta}$ and 
 $\lie{r}_{\,\R}=\lie{r}\cap\lie{k}\subseteq \lie{g}$.
The subspace $\lie{r}_{\,\R}$ defines a real structure on $\lie{r}$.

Let  $G$ be a connected and simply-connected complex Lie group with 
Lie algebra $\lie{g}$,
 $P$ the closed connected subgroup of $G$ with Lie algebra $\lie{p}$,
 % $X$ the homogeneous space $G/P$,
 $K$ the maximal compact subgroup of $G$ corresponding to $\lie{k}$,
 and $V$ the closed subgroup $K\cap P$ of $K$.
The coset space $X=G/P$ is a compact complex algebraic homogeneous space of
 positive Euler characteristic and every
 such space is of this form.
The inclusion of $K$ into $G$ induces a diffeomorphism from $K/V$ to
 $G/P$.
We assume in the following that $X=K/V$ is equipped with a $K$-invariant
 hermitian metric.
This metric induces the Hodge inner product on  
the space $A^\cdott(X)$ of smooth complex valued 
differential forms on $X$.
Let $\Delta=d^*d+dd^*$ be the associated Laplace operator on
$A^\cdott(X)$.
We equip the space of harmonic forms ${\mathcal H}^\cdott(X)=\ker(\Delta)$  
with the  induced hermitian metric from $A^\cdott(X)$.
The harmonic forms are contained in the subspace
 $A^\cdott(X)^K$ of $K$-invariant forms.
The natural inclusion from $A^\cdott(X)^K$ to $A^\cdott(X)$  induces 
a quasi-isomorphism of complexes
\[
\dis
(A^\cdott(X)^K,d)\longrightarrow(A^\cdott(X),d).
\]
The cohomology of these complexes can be calculated as follows.
The projection from $G$ to $X$ induces on (real) tangent spaces a surjection
 $T$ from $\wedge_\R\lie{g}$ to $\wedge_\R T_{e}X$. 
The restriction of $T$ defines an isomorphism 
between $\wedge_\R\lie{r}_{\,\R}$ and $ \wedge_\R T_{e}X$ \cite[Lemma 6.7]{2}.
Let $T_{e,\C}^*X=T_e^*X\otimes_\R\C$ be the space of all complex covectors
at the origin of $X$.
There is a unique isomorphism
\begin{equation}\label{tangent}
\dis
A:\wedge\lie{r}\rightarrow \wedge T_{e,\C}^*X
\end{equation}
defined so that $\left<Au,Tv\right>=(u,v)_\lie{g}$ 
holds for all $u\in \wedge \lie{r}$ and
$v\in \wedge_\R\lie{r}_{\,\R}$.
Here $(.,.)_\lie{g}$ denotes the bilinear form on $\wedge \lie{g}$ induced by the
Killing form of $\lie{g}$.

The Lie algebra $\lie{g}$ acts on $\wedge \lie{g}$ by the adjoint
representation.
The subspace $\wedge \lie{r}$ is stable under the restriction of this
representation to $\lie{g}_1$, i.e. $\wedge \lie{r}$ 
has the structure of a $\lie{g}_1$-module.
We denote the subspace $(\wedge \lie{r})^{\lie{g}_1}$ of
$\lie{g}_1$-invariant elements by $C$.
 %We observe that $C$ inherits a  hermitian metric from $\lie{r}$.
Let $d\in {\rm End}(\wedge \lie{g})$  be the coboundary operator on $\lie
{g}$, that is, the negative adjoint of the Chevalley-Eilenberg boundary 
operator on $\wedge \lie{g}$ with respect to the Killing form on $\lie{g}$.
The coboundary operator $d$ induces a differential $d$ on $C$.  
It is well known \cite[6.9]{2} that one obtains an isomorphism of
differential graded algebras 
\begin{equation}\label{400}
(A^\cdott(X)^K,d) \tilde{\longrightarrow}  (C,d)
\end{equation}
by mapping an invariant differential form $\omega$ to its restriction 
$\omega|_e\in (\wedge T_{e,\C}^*X)^{\lie{g}_1}=C$.
We obtain a canonical isomorphism of graded rings between $H^\cdott(X,\C)$ 
and $H(C,d)$.

The space $\lie{r}$ is not a Lie subalgebra of $\lie{g}$.
However the subalgebras 
$\lie{n}$ and $\lie{n}^{\theta}$ in the Lie algebra $\lie{g}$ 
equip $\lie{r}$ with a Lie algebra structure such 
that $[\lie{n},\lie{n}^{\theta}]=0$.
Let $\partialr\in{\rm End}(\wedge\lie{r})$ be the
Chevalley-Eilenberg boundary operator for the Lie algebra $\lie{r}$.
Let $b\in {\rm End}(\wedge\lie{r})$ be the corresponding coboundary
operator, that is, the negative adjoint  
of $\partialr$ with respect to the
restriction of $(.,.)_{\lie{g}}$ to $\wedge\lie{r}$.
The operators $b$ and $\partialr$ induce operators on $C$.
We consider the Laplacians
\begin{eqnarray*}
S&=&d\partialr+\partialr d\in {\rm End}(C)\\
L&=&b\partialr+\partialr b\in {\rm End}(C).
\end{eqnarray*}
In general, the operator $\partialr$ is not 
adjoint to $d$ with
respect to any hermitian metric on $X$ \cite[3.20]{EL}.
However Kostant shows in \cite[Section 4]{2} that $d$ 
and $\partialr$ are 
disjoint,
i.e. $d\partial(x)=0$ implies $\partial(x)=0$ and
$\partial d(y)=0$ implies $d(y)=0$ for all $x,y\in C$.
This implies that the kernel of $S$ computes the cohomologies $H(C,d)$
and $H(C,\partialr)$.
We obtain canonical isomorphisms
\[
\dis
\psi_{\Delta,S}:\ker(S)\tilde{\longrightarrow} H(C,d)
\tilde{\longrightarrow} \ker(\Delta)
\]
and
\[
\dis
\psi_{S,L}:\ker(L)\tilde{\longrightarrow} H(C,\partialr)
\tilde{\longrightarrow} \ker(S).
\]

Using the involution 
determined by the compact real form
 $\lie{k}$ of $\lie{g}$, we can define a positive definite hermitian 
inner product on $\wedge\lie{g}$ by
\[
\dis
\{u,v\}=(-1)^{\deg(u)}(u,v^{\theta})_\lie{g}
\]
for all $u,v\in \wedge \lie{g}$ \cite[3.3]{1}.
We equip the subspaces $\wedge{\lie{r}}$ and $C$ of $\wedge\lie{g}$ with the
induced hermitian inner product.
We are going to describe an orthogonal basis of 
the subspace $\ker(L)$ of $C$.
Therefore we consider the
representation of $\lie{g}_1$ on the Lie algebra homology
 $H_\cdott(\lie{n})$ induced by the adjoint action of $\lie{g}$ on
 $\wedge\lie{g}$.
This representation has the following description.
Let $\lie{h}$ denote the Cartan subalgebra $\lie{b}\cap\lie{b}^{\theta}$ of
$\lie{g}$, $R=R(\lie{g},\lie{h})$ the set of roots of $\lie{g}$
with respect to $\lie{h}$.
The choice of the Borel subalgebra $\lie{b}$ determines 
subsets $R_+$ and $R_-$ of $R$ 
of positive and negative roots respectively.
The maximal nilpotent ideal $\lie{n}$ is an $\lie{h}$-module under the
adjoint action of $\lie{h}$ on $\lie{g}$.
We denote by $R(\lie{n})$ the 
set of all roots whose root spaces lie in $\lie{n}$.
Let $W$ be the Weyl group of $\lie{g}$.
For every $\sigma\in W$, we have a subset
 $\Phi_\sigma=(\sigma R_-)\cap R_+$  of the set of roots $R$.
Corresponding to the parabolic subalgebra $\lie{p}$,
we define the set
\begin{equation}\label{w1}
W^1=\{\sigma\in W\ |\ \Phi_\sigma\subset R(\lie{n})\}.
\end{equation}
According to \cite[Corollary 8.1]{1} the $\lie{g}_1$-module 
$H_\cdott(\lie{n})$ may be decomposed as
\[
  H_\cdott(\lie{n})=\bigoplus_{\sigma\in
  W^1}M_\sigma
\]
where each $M_\sigma$ is an irreducible $\lie{g}_1$-module such that 
 $M_\sigma$ is not isomorphic to $M_\tau$ for $\sigma\neq \tau$.
The isomorphism $\wedge\lie{r}=\wedge\lie{n}\otimes\wedge\lie{n}^{\theta}$
induces an isomorphism \cite[Proposition 3.10]{2}
\[
\dis
H(C,\partialr)=(H_\cdott(\lie{n})\otimes H_\cdott(\lie{n}^{\theta}))^{\lie{g}_1}
\]
The Killing form of $\lie{g}$ puts $H_\cdott(\lie{n})$ and
$H_\cdott(\lie{n}^{\theta})$ in duality, so that $H_\cdott(\lie{n})$ is
the representation dual to $H_\cdott(\lie{n}^{\theta})$.
Using Schur's lemma, we get
\begin{eqnarray*}
(H_\cdott(\lie{n})\otimes H_\cdott(\lie{n}^{\theta}))^{\lie{g}_1}
&=&{\rm Hom}_{\lie{g}_1}(H_\cdott(\lie{n}),H_\cdott(\lie{n}))\\
&=&\oplus_{\sigma\in W^1}{\rm Hom}_{\lie{g}_1}(M_\sigma,M_\sigma)\\
&=&\oplus_{\sigma\in W^1}C_\sigma
\end{eqnarray*}
where each $C_\sigma$ is a one-dimensional space.
Let $h_\sigma$ be the preimage in $\ker (L)$ of a generator of $C_\sigma$ under
the isomorphism
\[
\dis
\ker(L)=H(C,\partialr)=(H_\cdott(\lie{n})\otimes
H_\cdott(\lie{n}^{\theta}))^{\lie{g}_1}=\bigoplus_{\sigma\in W^1}C_\sigma.
\]
The elements $h_\sigma$ 
form an orthogonal basis
of the subspace $\ker (L)$ of $C$ \cite[Theorem 5.4]{2}.
There is a canonical way to normalize the choice of $h_\sigma$
\cite[Proposition 5.5]{2} 
which is not needed in the following.
The class of a suitable 
multiple $\omega_\sigma$ of
 the image of $h_\sigma$ under the isomorphism
\begin{equation}\label{1400}
\ker(L)\stackrel{\psi_{S,L}}{\tilde{\longrightarrow}}\ker(S)
\stackrel{\psi_{\Delta,S}}{\tilde{\longrightarrow}}\ker(\Delta)
={\mathcal H}^\cdott(X).
\end{equation}
in  $H^\cdott(X,\C)$ is the Poincar\'e dual of the Schubert cell
in $X=G/P$ corresponding to $\sigma$ \cite[Theorem 6.15]{2}.
We call  $\om_{\sigma}$  the {\em Schubert form} corresponding to $\sigma$.
It follows from Schubert calculus that for each $\sigma\in W^1$ there is a
unique $\sigma'\in W^1$ such that 
\begin{equation}\label{dual}
\int_X\omega_\tau\wedge \omega_\sigma\neq 0\Leftrightarrow \tau=\sigma'
\end{equation}
holds for all $\tau\in W^1$.
We call $\omega_{\sigma'}$ the Poincar\'e dual form to $\omega_\sigma$
(an exact expression for $\sigma'$ is given in \cite[Corollary 2.6]{CN}).
Recall that the Hodge-$*$-operator on 
$A^\cdott(X)$ is determined by
\[
  \alpha\wedge*{\overline{\beta}}=\{\alpha,\beta\}{\mu_X},
\]
where $\mu_X=(n!)^{-1}\omega^n$ is the normalized top exterior power 
of the fundamental form $\omega$ of the metric.
Let us assume that the map \refg{1400} is an isometry.
Under this assumption, the Schubert forms $\omega_\sigma$, $\sigma\in W^1$,
form an orthogonal basis of ${\mathcal H}^\cdott(X)$.
Furthermore the Schubert forms are real and have the property \refg{dual}.
It follows that the Hodge $*$-operator maps a Schubert
form $\omega_\sigma$ to a multiple of $\omega_{\sigma'}$.

Let us finally assume that $X$ is an irreducible compact hermitian
symmetric space. 
We equip $X$ with its standard homogeneous hermitian metric. 
This is the unique $K$-invariant hermitian metric on $X$ for which 
\refg{tangent} becomes an isometry.
For compact hermitian symmetric spaces, this metric is K\"ahler and 
the space of harmonic forms  ${\mathcal H}^\cdott(X)$ coincides with 
the space  $A^\cdott(X)^K$ of $K$-invariant forms.
Furthermore, the Lie algebra $\lie{n}$ is commutative and 
the differentials $d$, $\partial$ and $b$ vanish on $C$.
We see in particular that $\psi_{S,L}$ is the identity 
and $\psi_{\Delta,S}$ coincides with \refg{400}.
It follows that \refg{1400} becomes
an isometry as \refg{400} is an isometry. 
Thus we have established:

\begin{prop}
\label{multprop}
Let $X$ be an irreducible
 hermitian symmetric space of compact type. Then the Hodge
star operator maps a Schubert form $\om_{\sigma}$ to a non-zero
multiple of the Poincar\'e dual form $\om_{\sigma'}$.
\end{prop}

The example in Section \ref{cexample} will show that the analogue of 
Proposition 
\ref{multprop} fails to hold for
arbitrary homogeneous spaces of the type $G/P$ considered above.

\section{Proof of the Main Theorem}
\label{proofsec}

Recall that for any
hermitian compact manifold  $X$ the 
Lefschetz operator $L:A^*(X)\ra A^{*+2}(X)$ on the space of smooth complex
valued differential forms is
given by $L(\eta)=\omega\wedge\eta$, where $\omega$ is the fundamental
form of the metric. When $X$ is a hermitian symmetric space $\omega$
is normalized to
coincide with the unique Schubert form of type $(1,1)$.

Let $\Lambda : A^{*+2}(X)\ra A^*(X)$ be the adjoint 
of $L$ with respect to the Hodge inner product on $A^*(X)$.
Recall that a differential form is {\em primitive} if
it lies in the kernel of $\L$. If $\dim_{\C}X=d$ and 
$\phi\in A^{p,p}(X)$ is a primitive form, then 
\begin{equation}
\label{wprim}
*\,L^r\phi=(-1)^p\frac{r!}{(d-2p-r)!}\,L^{d-2p-r}\phi.
\end{equation}
This follows from a more general theorem due to Weil
\cite[Chapter 1, Theorem 2]{Wei}. Applying (\ref{wprim})
when $\phi=1$ (and $p=0$) gives
\begin{equation}
\label{iterate}
*\,\om^r=\frac{r!}{(d-r)!}\,\om^{d-r}.
\end{equation}
The rest of the argument is a case by case analysis:

\medskip
\noin
(i) $X=G(m,n)$. In this case $\om=\Om_1$ and we have
\begin{equation}
\label{Adeg}
\Om_1^r=\sum_{|\l|=r} \frac{r!}{h_{\l}}\,\Om_{\l}
\end{equation}
(see for instance \cite[Example I.4.3]{M}). Using this in 
(\ref{iterate}) gives 
\begin{equation}
\label{Aequ}
\sum_{|\l|=r}\frac{*\,\Om_{\l}}{h_{\l}}=
\sum_{|\m|=mn-r}\frac{\Om_{\m}}{h_{\m}}.
\end{equation}
It follows from Proposition \ref{multprop} that 
$*\,\Om_{\l}=r_{\l}\Om_{\l'}$
for some real number $r_{\l}$. Since the $\Om_{\m}$ 
are linearly independent, (\ref{Aequ}) implies that $r_{\l}=
h_{\l}/h_{\l'}$, as required.

\medskip
\noin
(ii) $X=OG(n,2n)$. The analogue of equation (\ref{Adeg}) here is
\[
\dis
\Phi_1^r=\sum_{|\l|=r} \frac{r!}{g_{\l}}\,\Phi_{\l}.
\]
This follows from the Pieri formula for $X$ (see \cite{Hiller}, \cite{HB},
\cite[Section 6]{P}.) The rest of the argument is the same as in case (i).

\medskip
\noin
(iii) $X=LG(n,2n)$. The Pieri rule of \cite{HB} gives
\[
\dis
\Psi_1^r=\sum_{|\l|=r}
2^{\alpha(\l)}\frac{r!}{g_{\l}}\,\Psi_{\l}
\]
and we obtain the result as in the previous two cases.

\medskip
\noin
(iv) $X=Q(n)$. 
If $n=2k-1$ (resp.\ $n=2k$) then for $i\lequ k-1$ we have 
$\om^{n-i}
=2\,\om^{k-1-i}e$ (resp.\ $\om^{n-i}=2\,\om^{k-i}e_0=2\,\om^{k-i}e_1$) 
and the result follows immediately from (\ref{iterate}).
The formulas for $\omega^{n-i}$ are easily deduced 
from the formulas for $\omega^k$ in Section \ref{mt}. 
If $n=2k$ then
(\ref{iterate}) gives $*\,\om^k=\om^k$ and hence
\[
\dis
*\,e_0+*\,e_1=e_0+e_1.
\]
But Proposition \ref{multprop} implies that $*\,e_j$ is a multiple of
$e_{j+k}$ for $j=0,1$. Since $e_0$, $e_1$ freely generate $\H^{2k}(X)$,
we must have $*\,e_j=e_{j+k}$.
\endproof

\begin{figure}[htpb]
\epsfxsize 100mm
\center{\mbox{\epsfbox{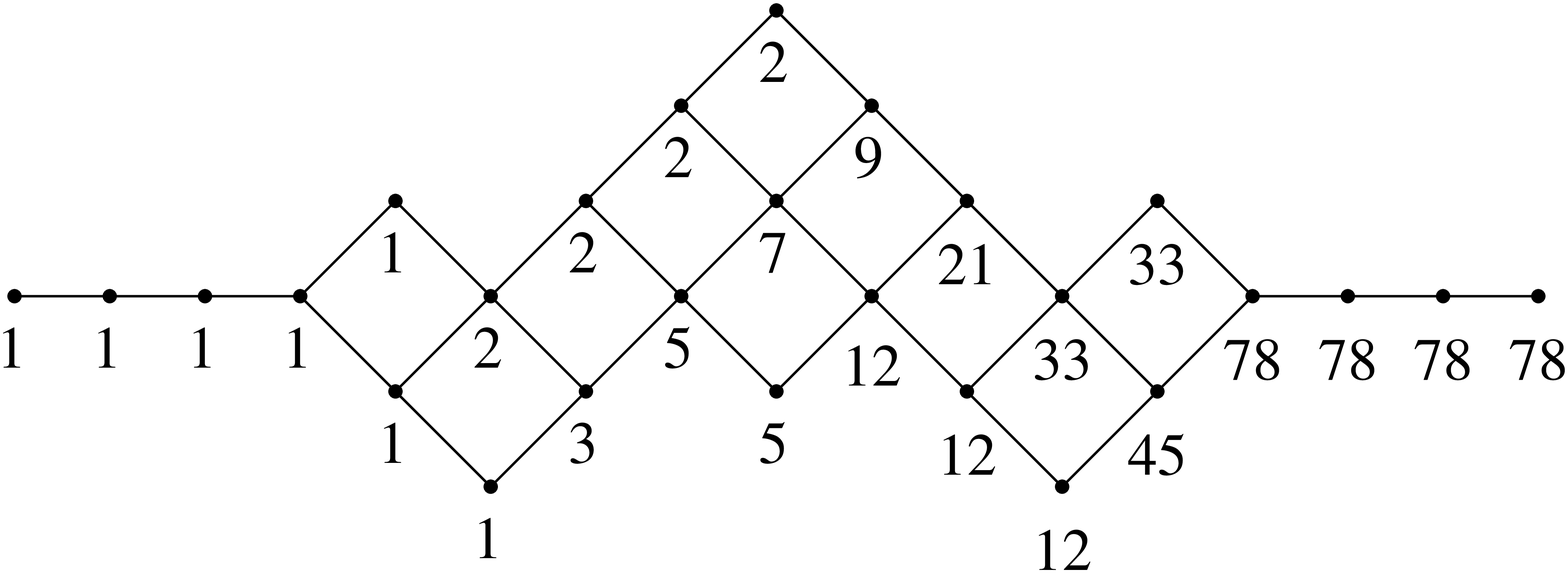}}} 
\caption{The poset of Schubert forms for $E_6/(SO(10)\cdot SO(2))$}
\label{e6}
\end{figure}

\medskip
The above argument applies to the exceptional spaces (v), (vi) of 
Section \ref{mt} as well. In general, the Schubert forms are parametrized 
by the Bruhat partially ordered set $W^1$, defined in (\ref{w1}). 
Each Schubert form
$\alpha$ corresponds to a node in the Bruhat poset; let $N(\alpha)$
denote the number of paths connecting $1$ to $\alpha$ in $W^1$. Figure
\ref{e6} shows the poset $W^1$ and the numbers $N(\alpha)$ in case
(v). If $\alpha$ is a $2|\alpha|$-form and $\alpha_1$ the unique
Schubert $(1,1)$-form then
\[
\dis
\alpha_1^r=\sum_{|\alpha|=r} N(\alpha) \alpha
\]
hence (\ref{iterate}) gives
\begin{equation}
\label{staralpha}
*\,\alpha = \frac{|\alpha|! \, N(\alpha')}{|\alpha'|! \, N(\alpha)} \,\alpha'
\end{equation}
where $\alpha'$ denotes the Poincar\'e dual of $\alpha$. Using
(\ref{staralpha}) one can compute the action of $*$ in examples (v)
and (vi).

\skipline
\noin
{\bf Remarks on normalization:} 1) (Forms)
In example (i) if we renormalize by setting $\wt{\Om}_{\l}=
\Om_{\l}/h_{\l}$ then Theorem \ref{mainthm} gives
$*\,\wt{\Om}_{\l}=\wt{\Om}_{\l'}$. However in this case 
$\wt{\Om}_{\l'}$ is not the Poincar\'e dual of $\wt{\Om}_{\l}$. 
A similar comment applies to the other spaces considered.

\medskip
\noin
2) (Metrics) Set $\om'=\rho\,\om$ for some $\rho>0$ and let 
$*$ (respectively $*'$) denote the Hodge star operator associated
with $\omega$ (respectively $\om'$). Then $*'=\rho^{d-k}\, *$
on the vector space $A^k(X)$ of differential $k$-forms on $X$. 
%This follows immediately from Weil's formula.

%\medskip
%\noin
%2) Wirtinger's Theorem (see \cite{GH}, p. 31) implies that
%the volume of a Schubert variety $X_{\l}$ in the Grassmannian 
%$G(m,n)$ satisfies
%\[
%\dis
%\mathrm{Vol}(X_{\l}) = \frac{1}{(\dim X_{\l})!}\deg(X_{\l}) = 
%\frac{1}{(mn-|\l|)!}\int_{X_{\l}}\Om_1^{mn-|\l|}=
%\frac{1}{h_{\l'}},
%\]
%and therefore
%\[
%\dis
%\mathrm{Vol}(X_{\l}) = \frac{h_{\l}}{h_{\l'}}\mathrm{Vol}(X_{\l'}).
%\]
%Since $*$ is an isometry, this shows that the volumes of a given 
%Schubert variety and its dual are
%in the same ratio as the norms squared of the corresponding Schubert forms.
%Similar remarks are true for the other hermitian symmetric spaces
%considered.

\section{Formulas for the Adjoint of the Lefschetz Operator}
\label{comb}
In this section we provide an explicit computation of the adjoint
$\L$ of the Lefschetz operator in the Grassmannian examples of
Section \ref{mt}. This calculation was done in a different way by Proctor
\cite{Pro1}, \cite{Pro2}, \cite{Pro3} in the case of minuscule
flag manifolds.
It follows from the definition of $\L$ that
\begin{equation}
\label{*L*}
\L=*L*
\end{equation}
as all non-zero harmonic forms occur in even degrees. Consequently one
can use our calculation for $*$ to find $\L$.

The Schubert forms in all our examples form a partially ordered set,
with the order induced from the Bruhat order on the underlying Weil
group (in combinatorial language they form an {\em irreducible 
Bruhat poset}, see
\cite[Section 2]{Pro2}.) The Lefschetz operator $L$ (respectively its
adjoint $\L$) is an order raising operator (resp. order lowering
operator) on the space of harmonic forms $\H^*(X)$. 
We will work out each of the three cases separately:

\medskip
\noin
(i) $X=G(m,n)$. The action of $L$ on the Schubert forms is given
by the Pieri rule:
\[
\dis
L(\Om_{\l})=\sum_{\m}\Om_{\m},
\]
the sum over all $\m\supset\l$ with $|\m|=|\l|+1$ (as usual, the
inclusion relation on partitions is defined by the containment
of diagrams.) Now equation (\ref{*L*}) and Theorem \ref{mainthm} give
\[
\dis
\L(\Om_{\l})=\sum_{\m}e_{\l\m}(m,n)\Om_{\m},
\]
the sum over all $\m$ with $\m\subset\l$ and $|\m|=|\l|-1$,
with $e_{\l\m}(m,n)=h_{\l}h_{\m'}/(h_{\l'}h_{\m})$.

\begin{prop}
\label{hookprop} For all $\m\subset\l$ with $|\m|=|\l|-1$, we have
\[
\dis
e_{\l\m}(m,n)=(m-i+\l_i)(n+i-\l_i),
\]
where $i$ is the unique index such that $\mu_i=\l_i-1$.
\end{prop}

\medskip
\noindent
{\bf Proof.} For any partition $\l$ with $\ell(\l)\lequ m$, the $\b$-{\em
sequence} $\b^{\l}$ is defined as the $m$-tuple
\[
\dis
\b^{\l}=(\l_1+m-1, \l_2+m-2, \ldots, \l_m+m-m).
\]
It is shown in \cite[Example I.1.1]{M} that
\[
\dis
h_{\l}=\frac{\prod_j\b^{\l}_j!}{\prod_{j<k}(\b^{\l}_j-\b^{\l}_k)}.
\]
Note that since
\[
\dis
\l'=(n-\l_m,n-\l_{m-1},\ldots,n-\l_1)
\]
we have
\[
\dis
\b^{\l'}=(m+n-1-\l_m,m+n-2-\l_{m-1},\ldots,m+n-m-\l_1).
\]
We now claim that
\[
\dis
\prod_{j<k}(\b_j^{\l}-\b_k^{\l})=\prod_{j<k}(\b_j^{\l'}-
\b_k^{\l'}).
\]
Indeed, the absolute value of these products is invariant under 
translation $t_p(x):=x+p$ and inversion $i(x):=-x$ of the entire
$\b$-sequence, and
\[
\dis
t_{m+n-1}( i(\b^{\l}))=\b^{\l'}.
\]
It follows that
\[
\dis
e_{\l\m}(m,n)=\frac{h_{\l}h_{\m'}}{h_{\l'}h_{\m}}=
\frac{\prod_j\b_j^{\l}!\prod_j\b_j^{\m'}!}
{\prod_j\b_j^{\l'}!\prod_j\b_j^{\m}!}.
\]
Finally, it is easy to check that
\[
\dis
\frac{\prod_j\b_j^{\l}!}{\prod_j\b_j^{\m}!}=m-i+\l_i \ \ \ \ 
\mathrm{and} \ \ \ \ 
\frac{\prod_j\b_j^{\m'}!}{\prod_j\b_j^{\l'}!}=
n+i-\l_i,
\]
as only one $\b$-number $(\b_i)$ changes (by one unit) when we pass
from $\l$ to $\m$. \endproof

\medskip
\noin
(ii) $X=OG(n,2n)$. In this case the Lefschetz operator satisfies
\[
\dis
L(\Phi_{\l})=\sum_{\m}\Phi_{\m},
\]
the sum over all strict partitions $\m$ with 
$\l\subset \m\subset \rho_{n-1}$ and $|\m|=|\l|+1$. Theorem
\ref{mainthm} and equation 
(\ref{*L*}) are now used to show that
%\begin{equation}
\[
\label{BLLpieri}
\L(\Phi_{\l})=\sum_{\m}f_{\l\m}(n)\Phi_{\m},
%\end{equation}
\]
the sum over $\m$ with $\m\subset\l$ and $|\m|=|\l|-1$,
where $f_{\l\m}(n)=g_{\l}g_{\m'}/(g_{\l'}g_{\m})$.

\begin{prop}
\label{hookprop2} For all strict $\m\subset\l$ with $|\m|=|\l|-1$, we have
\[
\dis
f_{\l\m}(n)=\left\{ \begin{array}{cl}
             n(n-1)/2  & \mbox{ if }  k = 1,\\
             n(n-1)-k(k-1) & \mbox{ otherwise,}
             \end{array} \right.
\]
where $k$ is the unique part of $\l$ which is not a part of $\m$.
\end{prop}

\medskip
\noindent
{\bf Proof.} The numbers $g_{\l}$ satisfy
\begin{equation}
\label{gequat}
g_{\l}=\prod_i\l_i!\cdot \frac{\prod_{i<j}(\l_i+\l_j)}
{\prod_{i<j}(\l_i-\l_j)}.
\end{equation}
This formula is due to Schur \cite{S}; see also \cite[Example III.8.12]{M}.
Now assume that $\l\smallsetminus\m = \{k\}$ and suppose 
that $k>1$. Let us compute the contribution of the three terms in 
(\ref{gequat}) to the quotient
$g_{\l}g_{\m'}/(g_{\l'}g_{\m})$:
the terms $\prod\l_i!$ contribute
\begin{equation}
\label{calc1}
\frac{\prod \l_i!\,\prod \m_i'!}{\prod \m_i!\,\prod \l_i'!}=
\frac{k!\,k!}{(k-1)!\,(k-1)!}=k^2.
\end{equation}
The contribution of the terms $\prod_{i<j}(\l_i+\l_j)$ is given by
\begin{equation}
\label{calc2}
\frac{\prod_{j\notin\{k,k-1\}}(k+j)}{\prod_{j\notin\{k,k-1\}}(k-1+j)}=
\frac{(k-1)(k+n-1)}{k^2}.
\end{equation}
The contribution of the terms $\prod_{i<j}(\l_i-\l_j)$ is given by
\begin{equation}
\label{calc3}
\frac{\prod_{j\neq k} |k-1-j|}{\prod_{j\neq k}|k-j|}=
\frac{n-k}{k-1}.
\end{equation}
Multiplying (\ref{calc1}), (\ref{calc2}) and (\ref{calc3}) together
gives
\[
\dis
f_{\l\m}(n)=\frac{g_{\l}g_{\m'}}{g_{\l'}g_{\m}}=
(n-k)(k+n-1).
\]
The case $k=1$ is handled similarly.
\endproof

\medskip
\noin
(iii) $X=LG(n,2n)$. The Lefschetz operator on $\H^*(X)$ is given by
\[
\dis
L(\Psi_{\l})=2\sum_{\m}\Psi_{\m} + \Psi_{\l^+}
\]
where the sum is over all (strict) $\m$ obtained from $\l$ by adding 
a box in a non-empty row and $\l^+=(\l_1,\ldots,\l_{\ell(\l)},1)$ (this
follows from the more general Pieri rule given by Hiller and Boe 
\cite{HB}). The computation of the adjoint operator $\L$ here is 
similar to the previous two cases, so we simply state the answer:
\[
\dis
\L(\Psi_{\l})=\frac{1}{2}\sum_{\m}f_{\l\m}(n+1)\Psi_{\m} +
f_{\l\l^-}(n+1)\Psi_{\l^-}
\]
where the sum is over all $\mu$ with $\ell(\m)=\ell(\l)$
obtained from $\l$ by subtracting a
box. The partition
$\l^-$ is defined to be empty if 1 is not a part of $\l$, and 
otherwise $\l^-=\l\smallsetminus 1$. 
%The numbers $f_{\l\m}$ are given by
%\[
%\dis
%f_{\l\m}=\left\{ \begin{array}{cl}
%             n(n+1)/2, & \mbox{ if }  k = 1\\
%             n(n+1)-k(k-1), & \mbox{ if }  k > 1
%             \end{array} \right.
%\]
%where $\l\smallsetminus\m = k$ defines $k$. 
Note that this calculation
was not given by Proctor, as $X$ is not a minuscule flag
manifold. 

We omit the computation of $\L$ for the remaining hermitian symmetric
spaces, which may be done in a similar fashion. The resulting
coefficients can be found in \cite{Pro3}.

\skipline
For any K\"ahler manifold $X$ define the
endomorphism $B:\H^*(X)\ra \H^*(X)$ by 
\[
\dis 
B=\sum_{i=0}^{2d}(d-i)
\mathrm{pr}_i,
\]
where $\mathrm{pr}_i$ is the projection onto the $i$th
homogeneous summand of $\H^*(X)$. It is well known that the operators
$L$, $\L$ and $B$ satisfy the commutator relations
%\begin{equation}
%\label{comrels}
\[
[B,L]=-2L, \ \qquad [B,\L]=2\L, \ \qquad  [\L,L]=B
%\end{equation}
\]
and hence we get a representation of ${\mathfrak s}
{\mathfrak l}(2,\C)$ on
$\H^*(X)$. It follows that $\L$ is completely determined by $L$
and $B$; see for instance \cite[Proposition 2]{Pro1} for a proof. 
Weil's formula (\ref{wprim}) and the Lefschetz decomposition theorem
now imply that $*:\H^*(X)\ra\H^*(X)$ 
is completely determined by the group
$\H^*(X)$ together with the action of the Lefschetz operator.
This allows us to include K\"ahler manifolds like example (ii$'$) in our
results.

%Finally, observe that one can use the above observation,
%Propositions \ref{multprop}--\ref{hookprop2} and 
%Proctor's computation of $\L$ to give an alternate
%proof of Theorem \ref{mainthm} which avoids 
%the use of (\ref{wprim}).

%
%$\clubsuit\clubsuit$
%Lemma 3.15 and Theorem 3.16 of \cite{Wel} should have interesting
%combinatorial consequences when applied to our situation. It should
%be possible to characterize the primitive forms in each example by
%using the relation
%\[
%\dis
%\Ker(\L:\H^k\ra \H^{k-2})=*\Ker(L:\H^{n-k}\ra\H^{n-k+2}).\ \ \ 
%\clubsuit\clubsuit
%\]

\section{An Example}
\label{cexample}

We calculate the action of the Hodge star operator for the complete
flag manifold
\[
F=F_{1,2,3}(\C^3)=SU(3)/S(U(1)^3)
\] 
which
parametrizes complete flags in a three dimensional complex vector space.
We will see that the analogue of Proposition \ref{multprop} fails when
$F$ is equipped with any $SU(3)$-invariant metric.

There is a universal vector bundle $E$ over $F$ together with a 
tautological filtration
\[
0=E_0\subsetneq E_1 \subsetneq E_2 \subsetneq E_3=E
\]
by subbundles such that each quotient $L_i=E_i/E_{i-1}$ is a line bundle on
$F$.
We consider the natural $\C$-algebra homomorphism
from $\C\,[x_1,x_2,x_3]$ to
 $H^*(F,\C)$ which maps $x_i$ to $y_i=-c_1(L_i)\in H^2(F,\C)$.
It is well known that this map is surjective and that its kernel is
generated by the elementary symmetric polynomials
 $x_1+x_2+x_3$, $x_1x_2+x_1x_3+x_2x_3$, and $x_1x_2x_3$.
This yields relations
 $y_3=-y_1-y_2$, $y_2^2=-y_1^2-y_1y_2$, and $y_1^3=0$
in $H^*(F,\C)$.
The cycle classes of the Schubert varieties in $F$ define an $\C$-basis of
 $H^*(F,\C)$.
This basis is given by the classes $1$, $y_1$, $y_1+y_2$, $y_1^2$,
 $y_1y_2$, and $y_1^2y_2$ (these are the {\em Schubert polynomials} for
$S_3$; see \cite{LS}).

%The standard Riemannian metric on $F$ is the normal
%metric \cite[7.86, 8.13]{Be} induced by the opposite of the Killing form on
%$\lie{s}\lie{u}(3)$; $F$ carries a unique homogeneous 
%hermitian metric whose
%underlying Riemannian metric is the standard one.

Equip $F$ with any $SU(3)$-invariant hermitian metric and 
denote the fundamental form of this metric by $\omega$.
We will show that the Hodge-$*$-operator associated
with $\omega$ satisfies
%\begin{equation}\label{example}
\[
*\,y_1=\lambda\, y_1^2+\mu\,(y_1y_2)
%\end{equation}
\]
for some $\lambda\neq 0$.
Observe that this is in contrast with our result in the hermitian
symmetric space case as $y_1y_2$ is Poincar\'e dual to $y_1$. 

We will work with the differential $(1,1)$-forms 
$\Omega_{ij}$, $1\lequ i< j\lequ  3$ on $F$
constructed in \cite[Section 5]{Ta}.
The $\Om_{ij}$ are a basis of the $3$-dimensional
space of $SU(3)$-invariant forms in $A^2(F)$; note that there are no
$SU(3)$-invariant 1-forms on $F$.
The subspace of harmonic forms ${\mathcal H}^2(F)$ in $A^{2}(F)$ has 
dimension two. We have 
\[
\dis
\omega=\a\,\Om_{12}+ \b\, \Om_{13}+\g\, \Om_{23}
\]
for some positive real numbers $\a$, $\b$, $\g$; 
the metric is K\"ahler if and only if $\b=\a+\g$ (see
for instance \cite[Section 4]{Bo}).

%The validity of equation \refg{example} for some $\lambda\neq 0$ is not 
%affected by a scaling of the metric.

The harmonic representatives $h_i$ of the classes
$y_i$ do not depend on the choice of invariant metric, and are given as
\cite[Corollary 3]{Ta}
\[
h_1=\Omega_{12}+\Omega_{13}\,\,,\,\,
h_2=-\Omega_{12}+\Omega_{23}\,\,,\,\,
h_3=-\Omega_{13}-\Omega_{23}.
 \]
One checks easily that $z=\a\,\Om_{12}-\b\,\Om_{13}$ 
and $z'=\g\,\Omega_{23}-\a\,\Omega_{12}$
are $\omega$-primitive, i.e. satisfy $\omega^2\wedge z=0$ and
$\omega^2\wedge z'=0$ in $A^\cdott(F)$.
Using formulas \refg{wprim}, \refg{iterate}, and the equalities
\[
3\,\a\,\Omega_{12}=\omega+z-z'\,\,,\,\,
3\,\b\,\Omega_{13}=\omega-2\,z-z'\,\,,\,\,
3\,\g\,\Omega_{23}=\omega+z+2\,z'\,,
\]
one calculates
\[
*\,\a\,\Omega_{12}=\b\g\,\Omega_{13}\wedge\Omega_{23}\,\,,\,\,
*\,\b\,\Omega_{13}=\a\g\,\Omega_{12}\wedge\Omega_{23}\,\,,\,\,
*\,\g\,\Omega_{23}=\a\b\,\Omega_{12}\wedge\Omega_{13}
\]
and
\[
*\,h_1=\frac{\a\g}{\b}\,\Omega_{12}\wedge \Omega_{23}+
\frac{\b\g}{\a}\,\Omega_{13}\wedge \Omega_{23}.
\]
We conclude that 
\[
\l=\int_F(*h_1)\wedge (h_1+h_2)=\frac{\a\g}{\b}
\int_F \Omega_{12}\wedge\Omega_{13}\wedge
\Omega_{23}\neq 0.
\]

\end{document}